\documentclass[]{amsart}

\title{Non-productive duality properties of topological groups}
\author{Masasi Higasikawa}
\address{Eda Laboratory, School of Science and Engineering, Waseda University,
Shinjuku-ku, Tokyo, 169-8555, Japan}
\email{higasik@logic.info.waseda.ac.jp}
\keywords{Topological group, character, duality, sequential topology,
exponential Diophantine equation.}
\subjclass[2000]{Primary 22A05, 43A40; Secondary 11D61, 11Z05} 
\date{October 26, 2000}
\thanks{Work partly supported by JSPS Research Fellowships for
Young Scientists and a Grant-in-Aid for Scientific Research
from the Ministry of Education, Science and Culture of Japan.}

\theoremstyle{plain}
\newtheorem{theorem}{Theorem}[section]
\newtheorem{prop}[theorem]{Proposition}
\newtheorem{lemma}[theorem]{Lemma}

\theoremstyle{remark}
\newtheorem{remark}[theorem]{Remark}

\newcommand{\Pair}[1]{\left\langle #1\right\rangle}
\newcommand{\N}{{\mathbf N}}
\newcommand{\Z}{{\mathbf Z}}
\newcommand{\R}{{\mathbf R}}
\newcommand{\T}{{\mathbf T}}
\newcommand{\Num}{\in\N}
\newcommand{\Id}{\mathrm{id}}

\newcommand{\XI}{{\rm X(1)}}
\newcommand{\XII}{{\rm X(2)}}

\newcommand{\Em}{}

\begin{document}

\begin{abstract}
We address two properties for Abelian topological groups: ``every closed
subgroup is dually closed'' and ``every closed subgroup is dually embedded.''
We exhibit a pair of topological groups such that each has both of the
properties and the product has neither, which refutes a remark of
N.~Noble. 

These examples are the additive group of integers topologized with respect to
a convergent sequence as investigated by E.G.~Zelenyuk and I.V.~Protasov. The
proof for the product relies on a theorem on exponential Diophantine
equations.
\end{abstract}

\maketitle

\section{Introduction}
Several duality properties are well-established for locally-compact Abelian
groups. Although some are shared beyond the class (\cite{BHM}, \cite{Ka},
\cite{No}, \cite{Va}), the boundaries may be rather vague. Our concern is for
two of them: ``every closed subgroup is dually closed'' and ``every closed
subgroup is dually embedded.'' We denote the former by \XI\ and the latter by
\XII\ after \cite{BHM}.

Some definitions and conventions follow. All topological groups treated here
are Hausdorff and Abelian, and a character is a continuous homomorphism into
the torus $\T=\R/\Z$, unless otherwise stated. A subgroup $H$ of a topological
group $G$ is \Em{dually closed} if for each $g\in G\setminus H$, there exists
a character $\chi$ of $G$ that separates $H$ and $g$, i.e., $\chi$ is 
identically zero on $H$ and $\chi(g)\neq 0$. We say that $H$ is
\Em{dually embedded} if each charater of $H$ extends to one of $G$.

That a locally-compact group has \XI\ and \XII\ is part of the celebrated
Pontryagin-van Kampen duality theorem. S.~Kaplan \cite{Ka} shows the same for
a product of locally-compact groups. Erroneously referring to the proofs,
N.~Noble claims in Section 3 of \cite{No} that \XI\ and \XII\  are each
preserved under arbitrary products.

\begin{remark}
Kaplan's arguments only show that a product $\prod_{i\in J}G_i$ of topological
groups has \XI\ (\XII, respectively) \Em{whenever} each subproduct
$\prod_{i\in J}G_i$ does for finite $J\subseteq I$. So the problem is reduced
to the preservation under finite products, which is out of question for
locally-compact groups.
\end{remark}

We exhibit a counterexample to Noble's assertion above. Let $\Z\{2^n\}$ denote
the additive group of integers with the strongest group topology such that
$2^n$ converges to $0$, and $\Z\{3^n\}$ similarly (E.G.~Zelenyuk and
I.V.~Protasov \cite{ZP}). Each has both of \XI\ and \XII\ (Proposition
\ref{factor}). We show that the product $\Z\{2^n\}\times\Z\{3^n\}$ has neither
(Theorem \ref{conc}) as to the diagonal. The proof is dependent on the fact
that the diagonal is closed and discrete (Theorem \ref{diag}), which we get by
invoking a theorem on exponential Diophantine equations (cf. \cite{Sc}).

In Section 2, we investigate some basic properties of the topological group
$\Z\{p^n\}$. Section 3 is for number theory necessary for the rest of the
argument. The proof of the non-preservation of \XI\ and \XII\ is completed in
Section 4. In Section 5, we construct another counterexample with metrics which
have close relation to certain family of exponential Diophantine equations.

\section{Group topologies with a sequence convergent}

Let $G$ be a group, Abelian or not, and $\Pair{a_n:n\Num}$ a sequence in
$G$. Then there exists the strongest group topology on $G$ such that
$\Pair{a_n:n\Num}$ converges to the neutral element. We denote by $G\{a_n\}$
the topological group $G$ with this topology, which need not be Hausdorff.
Zelenyuk and Protasov \cite{ZP} investigates such topological groups
determined by a convergent filter in Abelian and Hausdorff case; unfortunately
their recently published monograph \cite{PZ} have not been available for the
author.

The strongest topology thus defined has a useful characterization as follows.

\begin{prop}\label{cont}
Suppose that $G$ and $\Pair{a_n:n\Num}$ are as above. Let $H$ be a topological
group and $f:G\rightarrow H$ a homomorphism between abstract groups. Then $f$
is continuous if and only if $f(a_n)\rightarrow 1$ in $H$.
\qed
\end{prop}

Our main concern is for topological groups of the form $\Z\{p^n\}$ with $p$
prime. Their closed subgroups and characters are explicitly described in
\cite{ZP}.

\begin{prop}\label{factor}
\begin{enumerate}
\item The closed subgroups of $\Z\{p^n\}$ are $p^m\Z$ ($m=0,1,...$) and
 $\{0\}$.
\item The character group of $\Z\{p^n\}$ is (identified with)
 $\Z[1/p]/\R\subset\T$, where $\Z[1/p]$ is the ring generated by $1/p$ over
 $\Z$.
\item Both of \XI\ and \XII\ hold for $\Z\{p^n\}$.
\qed
\end{enumerate}
\end{prop}

Let $p$ and $q$ be distinct primes. The characters of the product
$\Z\{p^n\}\times\Z\{q^n\}$ are induced by those of the factors, and hence we
see how they behave to the diagonal $\Delta=\{\Pair{u,u}:u\in\Z\}$.

\begin{lemma}\label{prodchar}
\begin{enumerate}
\item The diagonal and each element lying outside cannot be separated by the
 characters.
\item The group of characters of $\Delta$ extendable to the whole product is
 a proper subgroup of $\T$.
\end{enumerate}
\end{lemma}
\begin{proof}
Each character $\chi\in\Z\{p^n\}\times\Z\{q^n\}$ is of form
$\chi(u,v)=\chi(1,0)u+\chi(0,1)v$ with $\chi(1,0)\in\Z[1/p]/\R$ and
$\chi(0,1)\in\Z[1/q]/\R$.
\end{proof}

Therefore \Em{assuming} that the diagonal is closed and discrete
(Theorem \ref{diag}), we have the desired non-preservation result.

\begin{theorem}\label{conc}
The product $\Z\{p^n\}\times\Z\{q^n\}$ has neither \XI\ nor \XII.
\qed
\end{theorem}

\begin{remark}
By \cite[Theorem 6]{ZP}, $\Z\{p^n\}$ is not Fr\'echet-Urysohn, so not
metrizable a fortiori. We observe that it is not $\alpha_4$ (cf. \cite{A1},
\cite{A2}) due to \cite[Lemma 3]{ZP}; Nyikos \cite{Ny} proves that
Fr\'echet-Urysohn group is $\alpha_4$.
\end{remark}

\section{$S$-unit equations}

We recall a finiteness theorem for exponential Diophantine equations, a
special case of \cite[Ch.~V, Theorem 2A]{Sc}. Let $S$ be a finite set of
primes.
A rational number is said to be an $S$-unit if it belongs to the multiplicative
group generated by $S\cup\{-1\}$. The set of $S$-units is denoted by $U_S$.
\begin{theorem}\label{sunit}
Up to scalar multiplications, the equation $x_1+\cdots+x_k=0$ has only finitely
many solutions $\Pair{x_1,...,x_k}$ in $S$-units whose non-trivial subsums do
not vanish.
\qed
\end{theorem}

As a corollary, we have finiteness for subsums as to a kind of $S$-unit
equation, which we need in the next section.

\begin{lemma}\label{sunitc}
Suppose that $S$ and $T$ are disjoint finite sets of primes. Let the tuple
$\Pair{x_1,...,x_k,y_1,...,y_l}$ runs through the solutions of the equation
$x_1+\cdots+x_k=y_1+\cdots+y_l$ with $x_i\in U_S\cup\{0\}$ and
$y_j\in U_T\cup\{0\}$ for $1\leq i\leq k,1\leq j\leq l$. Then the sum
$x_1+\cdots+x_k$ has only finitely many values.
\end{lemma}
\begin{proof}
For notational convenience, we consider the equation
$x_1+\cdots+x_k+x_{k+1}+\cdots+x_{k+l}=0$ replacing each $y_j$ with $-x_{k+j}$.
For $I\subseteq\{i:1\leq i\leq k+l\}$, we say a tuple $\Pair{x_i:i\in I}$ to be
an $I$-solution and $\sum_{1\leq i\leq k,i\in I}x_i$ an $I$-admisible sum if
$\sum_{i\in I}x_i=0$ with $\{x_i:1\leq i\leq k,i\in I\}\subseteq U_S\cup\{0\}$
and $\{x_i:k+1\leq i\leq k+l,i\in I\}\subseteq U_T\cup\{0\}$. We shall show
that there are only finitely many $I$-admissible sums.

Note that if $I$ is contained in $\{i:1\leq i\leq k\}$ or
$\{i:k+1\leq i\leq k+l\}$, which we call a degenerate case, then the only
$I$-admissible sum is $0$.

We show the finiteness by the induction on the size of $I$. If $|I|=1$, then
it is degenerate.

Assume $|I|\geq 2$ and it is not degenerate. Then the equivalence among
$I$-solutions up to scalar multiplications is just the identity. An
$I$-admissible sum is either as to an $I$-solution whose non-trivial subsums
do not vanish or of the form $u+u'$ such that $u$ is a $J$-admissible sum and
$u'$ is an $I\setminus J$-admissible sum for a non-trivial subset $J$ of $I$.
By the theorem and by the induction hypothesis respectively, there are only
finitely many such sums.
\end{proof}

\section{Reduction to number theory}

We investigate the topological properties of $\Z\{p^n\}\times\Z\{q^n\}$ so
that the result in the previous section may be applicable.

\begin{lemma}\label{prodseq}
For topological groups of the form $G\{a_n\},H\{b_n\}$, which need not be
Abelian or Hausdorff, the product topology $G\{a_n\}\times H\{b_n\}$ is the
same as $(G\times H)\{c_n\}$, where
$c_{2n}=\Pair{a_n,1},c_{2n+1}=\Pair{1,b_n}$.
\end{lemma}
\begin{proof}
It is straightforward that $c_n\rightarrow\Pair{1,1}$ in
$G\{a_n\}\times H\{b_n\}$. So the product has at most as strong topology as
that determined by $\Pair{c_n:n\Num}$. Hence it is sufficient to show that the
identity map $\Id:G\{a_n\}\times H\{b_n\}\rightarrow(G\times H)\{c_n\}$ is
continuous.

In general, the product topology of topological groups $G$ and $H$ has
following characterization. Let $i:G\rightarrow G\times H$ and
$j:H\rightarrow G\times H$ be natural injections:
$i(g)=\Pair{g,1},j(h)=\Pair{1,h}$. Suppose that $K$ is a topological group and
$f:G\times H\rightarrow K$ is a homomorphism between abstract groups. Then $f$
is continuous if and only if so are both of $f\circ i$ and $f\circ j$.

Due to Proposition \ref{cont},
$\Id\circ i:G\{a_n\}\rightarrow(G\times H)\{c_n\}$ and
$\Id\circ j:H\{b_n\}\rightarrow(G\times H)\{c_n\}$ are continuous. So we are
done.
\end{proof}

We depend on the following two results in \cite{ZP} concerning Hausdorff
abelian topologial groups of the form $G\{a_n\}$.

\begin{lemma}[{\cite[Lemma 2]{ZP}}]\label{conv}
If $g_m\rightarrow 0$ in $G\{a_n\}$, then there exists a positive integer $k$
such that
$g_m\in\{x_1+\cdots+x_k:(\forall i)(x_i\in\{\pm a_n:n\Num\}\cup\{0\})\}$ for
sufficiently large $m$.
\qed
\end{lemma}

\begin{theorem}[{\cite[Theorem 7]{ZP}}]\label{seq}
The topology of $G\{a_n\}$ is sequential.
\qed
\end{theorem}

Now the last piece of the proof follows. Note that every discrete subgroup of
a Hausdorff topological group is closed.

\begin{theorem}\label{diag}
The diagonal $\Delta$ is discrete (and closed) in $\Z\{p^n\}\times\Z\{q^n\}$.
\end{theorem}
\begin{proof}
Applying Theorem \ref{seq} to $\Z\{p^n\}\times\Z\{q^n\}$
(Lemma \ref{prodseq}), we observe that the product is also sequential.
So it is sufficient to show that $\Delta$ has no nontrivial convergent
sequence.

Suppose that $\Pair{\Pair{u_i,u_i}:i\Num}$ is a sequence in $\Delta$
converging to $\Pair{0,0}$. By Lemma \ref{conv}, there exists $k$ such that
$u_i$ is a sum of less than $k$ numbers in $\{\pm 2^n:n\Num\}$ and in
$\{\pm 3^n:n\Num\}$, respsctively, for sufficiently large $i$. Due to
Lemma \ref{sunitc}, there are only finitely many such sums. Therefore
$u_i$ is eventually equal to $0$.
\end{proof}

\begin{remark}
Sequentiality need not be preserved under direct products even for
topological groups (\cite{To}, \cite{Sh}).
\end{remark}

\section{Metric counterexample}
Proposition \ref{factor} holds not only for $\Z\{p^n\}$ but also for a group
topology such that each $p^m\Z$ is closed and $p^n\rightarrow 0$. Among those
topologies the strongest is that of $\Z\{p^n\}$ and the weakest is the $p$-adic
topology. They have the same closed (open) subgroups and characters.

If $\mathcal T$ is such a topology and $\mathcal U$ is one for another prime
$q\neq p$, then Lemma \ref{prodchar} is true for
$\Pair{\Z,\mathcal T}\times\Pair{\Z,\mathcal U}$ as well. So we have another
counterexample, which we shall construct in this section, provided that the
diagonal is discrete.

Suppose that $\delta:\{p^n:n\Num\}\rightarrow\R_{>0}$ is
a non-increasing function with $\delta(p^n)\rightarrow 0$. We define a function
$||\cdot||_\delta:\Z\rightarrow\R$ by
\[||u||_\delta=\inf\left\{\sum_i \delta(p^{n_i}):u=\sum_i e_i p^{n_i},
 e_i\in\{1,-1\}, n_i\in\N\right\}.\]
We denote by $\Z_\delta$ the topological group $\Z$ with the metric induced by
$||\cdot||_\delta$. This topology is strictly between the $p$-adic topology and
the strongest topology.

For distinct primes $p,q$, we construct $\delta:\{p^n:n\Num\}\rightarrow\R$
and $\epsilon:\{q^n:n\Num\}\rightarrow\R$ with certain property. We need
another corollary to Theorem \ref{sunit}.

\begin{lemma}
We consider the equation $x_1+\cdots+x_k=y_1+\cdots+y_l$ in Lemma \ref{sunitc}
under the restriction such that no (non-empty) subsum of
$x_1+\cdots+x_k$ or of $y_1+\cdots+y_l$ vanishes, in particular
$x_i\in U_S,y_j\in U_T$. Then the number of solutions is finite.
\end{lemma}
\begin{proof}
We proceed similarly as in the proof of Lemma \ref{sunitc}. We change the
definition of $I$-solutions according to the restriction: they are tuples
$\Pair{x_i:i\in I}$ such that $\sum_{i\in I}x_i=0$,
$\{x_i:1\leq i\leq k,i\in I\}\subseteq U_S$,
$\{x_i:k+1\leq i\leq k+l,i\in I\}\subseteq U_T$ and no subsum of
$\sum_{1\leq i\leq k,i\in I}x_i$ or of $\sum_{k+1\leq i\leq k+l,i\in I}x_i$
vanishes. So there is no $I$-solutions in the degenerate case.

By the induction on the size of $I$, we prove the finiteness of $I$-solutions.
Assume that $I$ is not degenerate. Then there are two possibilities for an
$I$-solution: either it has no vanishing non-trivial sums or it is the union
of a $J$-solution and an $I\setminus J$-solution for a non-trivial subset $J$
of $I$. Hence Theorem \ref{sunit} and the induction hypothesis yields the
conclusion.
\end{proof}

Suppose that a non-zero integer $u$ has representations as
$u=\sum_{1\leq i\leq k} e_i p^{m_i}=\sum_{1\leq j\leq l} f_j q^{n_j}$ with
$e_i,f_j\in\{1,-1\}$ and $m_i,n_j\in\N$. To estimate $||u||_\delta$ and
$||u||_\epsilon$, we may assume that no subsum is zero. Due to the lemma above,
for fixed $k,l$, there are only finitely many summands $p^{m_i}$ which appear
in such representations for some non-zero integer $u$. Let  $F(p,q,k,l)$ denote
the finite set of such $p^{m_i}$. For a  positive integer $s$, set
$F(p,q,s)=\bigcup_{k,l\leq s}F(p,q,k,l)$. Then it is a finite set
non-decreasing with respect to $s$ and $\bigcup_{s\Num}F(p,q,s)=\{p^n:n\Num\}$.

Now we define particular $\delta$ and $\epsilon$ by
\[\delta(p^n)=1/\min\{s:p^n\leq\max F(p,q,s)\},\]
\[\epsilon(q^n)=1/\min\{s:q^n\leq\max F(q,p,s)\}.\]
Then they are non-increasing and tend to $0$. Moreover we have
\[p^n\in F(p,q,s)\Rightarrow \delta(p^n)\geq 1/s,\]
\[q^n\in F(q,p,s)\Rightarrow \epsilon(q^n)\geq 1/s.\]

We estimate $||u||_\delta$ and $||u||_\epsilon$ for a integer $u$ with the
representations above. If we set $s=\max\{k,l\}$, then
$\delta(p^{m_i}),\epsilon(q^{n_j})\geq 1/s$ for each summand. Without loss of
generality, we may assume $k\geq l$. Then we have
\[\sum_{1\leq i\leq k}\delta(p^{m_i}) \geq k\cdot(1/s) =1,\]
and hence $||u||_\delta \geq 1$.

Accordingly it does not occur that both $||u||_\delta$ and $||u||_\epsilon$ are
less than $1$ for a non-zero integer $u$. Therefore the diagonal is discrete
in $\Z_\delta\times\Z_\epsilon$. 

\begin{theorem}
Neither \XI\ nor \XII\ is not preserved under the product
$\Z_\delta\times\Z_\epsilon$ for $\delta$ and $\epsilon$ decreasing slowly
enough.
\qed
\end{theorem}

\end{document}